\documentclass[11pt]{article}
\usepackage{amssymb}

\usepackage{graphicx}
\usepackage{amsmath}
\usepackage{makeidx}
\usepackage{indentfirst}

\newcounter{resultnum}[section]
\newtheorem{conclusion}{Conclusion}[section]

\newcounter{conclusionnum}[section]

\newcounter{conditionnum}[section]

\newcounter{conjecturenum}[section]

\newcounter{examplenum}[section]

\newcounter{exercisenum}[section]
\newtheorem{lemma}{Lemma}[section]

\newcounter{lemmanum}[section]

\newcounter{notationnum}[section]
\newtheorem{theorem}{Theorem}[section]

\newcounter{theoremnum}[section]
\newtheorem{definition}{Definition}[section]

\newcounter{definitionnum}[section]
\newtheorem{corollary}{Corollary}[section]

\newcounter{corollarynum}[section]
\newtheorem{remark}{Remark}[section]

\newcounter{remarknum}[section]
\newtheorem{proposition}{Proposition}[section]

\newcounter{propositionnum}[section]

\newcounter{acknowledgementnum}[section]

\newcounter{algorithmnum}[section]

\newcounter{axiomnum}[section]

\newcounter{casenum}[section]
\newtheorem{claim}{Claim}[section]

\newcounter{claimnum}[section]

\newcounter{summarynum}[section]

\newcounter{problemnum}[section]
\newenvironment{proof}[1][]{\textbf{Proof.} }{}

\begin{document}

\title{The Entropy of Lagrange--Finsler \\
Spaces and Ricci Flows}
\date{October 18, 2008}
\author{ Sergiu I. Vacaru\thanks{
 Sergiu.Vacaru@gmail.com } \\
%EndAName
{\quad} \\
\textsl{The Fields Institute for Research in Mathematical Science} \\
\textsl{222 College str., 2d Floor, Toronto, ON, Canada M5T 3J1} \\
and \\
\textsl{Faculty of Mathematics, University "Al. I. Cuza" Ia\c si}, \\
\textsl{\ 700506, Ia\c si, Romania} }
\maketitle

\begin{abstract}
We formulate a statistical analogy of regular Lagrange mechanics and Finsler
geometry derived from Grisha Perelman's functionals and generalized for
nonholonomic Ricci flows.  Explicit constructions are elaborated when
nonholonomically constrained flows of Riemann metrics result in Finsler like
configurations, and inversely, when geometric mechanics is modelled on
Riemann spaces with a preferred nonholonomic frame structure.

\vskip0.3cm



\end{abstract}

%\tableofcontents

\section{Introduction}

The Ricci flow theory became a very powerful method in understanding the
geometry and topology of Riemannian manifolds \cite{ham1,gper1,gper2,gper3}
(see also reviews \cite{caozhu,kleiner,rbook} on Hamilton--Perelman theory
of Ricci flows). There were proposed a number of important innovations in
modern physics and mechanics.

Any regular Lagrange mechanics and analogous gravity theory can be naturally
geometrized on nonholonomic Riemann manifolds as models of Lagrange, or
Finsler, spaces \cite{vnhrf1,vijmpa1}, see Refs. \cite{ma1,ma2,bejf,vsgg}
for details and applications to modern physics. One of the major goals of
geometric mechanics is the study of symmetry of physical systems and its
consequences. In this sense, the ideas and formalism elaborated in the Ricci
flow theory provide new alternatives for definition of 'optimal' geometric
configurations and physical interactions.

A Riemannian geometry is defined completely on a manifold\footnote{%
for simplicity, in this work we shall consider only smooth and orientable
manifolds} provided with a symmetric metric tensor and (uniquely defined to be
metric compatible and torsionless) Levi--Civita connection structures.
Contrary, the Lagrange and Finsler geometries and their generalizations are
constructed from three fundamental and (in general) independent geometric
objects: the nonlinear connection, metric and linear connection. Such models
were developed when the main geometric structures are derived
canonically from a fundamental effective, or explicit, Lagrange (Finsler)
function and have an alternative realization as a Riemann geometry with
a preferred nonholonomic frame structure. Following such ideas, in Ref. \cite%
{vnhrf1}, we proved that Ricci flows of Riemannian metrics subjected to
nonholonomic constraints may result in effective Finsler like geometries and
that any Lagrange--Finsler configuration can be 'extracted' from the corresponding
nonholonomic deformations of frame structures. An important result is that
the G. Perelman's functional approach \cite{gper1,gper2,gper3} to Ricci
flows can be redefined for a large class of canonical metric compatible
nonlinear and linear connections. For regular Lagrange systems, this allows
us not only to\ derive the evolution equations and establish certain optimal
geometric and topological configurations but also to construct canonical
statistical and thermodynamical models related to effective mechanical,
gravitational or gauge interactions.

The aim of this paper is to analyze the possible applications of the theory of
Ricci flows to geometric mechanics and related thermodynamical models. We
shall follow the methods elaborated in Sections 1-5 of Ref. \cite{gper1}
generalizing this approach to certain classes of Lagrange and Finsler metrics
and connections (see a recent review on the geometry of nonholonomic
manifolds and locally anisotropic spaces in Ref. \cite{vrflg,vsgg}). It
should be emphasized that such constructions present not only a geometric
extension from the canonical Riemannian spaces to more sophisticate
geometries with local anisotropy but launch a new research program \cite%
{vnhrf1,vijmpa1,vv1,vv2} on Ricci flows of geometric and physical objects
subjected to nonholonomic constraints.

The paper is organized as follows: In section 2, we outline the main results
on metric compatible models of Lagrange and Finsler geometry on nonholonomic
manifolds. \ The G. Perelman's functional approach to Ricci flow theory is
generalized for Lagrange and Finsler spaces in section 3. \ We derive the
evolution equations for Lagrange--Ricci systems in section 4. A statistical
approach to Lagrange--Finsler spaces and Ricci flows is proposed in section
5. Finally, we discuss the results in the section 6. Some relevant
formulae are presented in the Appendix.

\section{Lagrange Mechanics and \newline
N--anholonomic Manifolds}

Let us consider a manifold $\mathbf{V,}$ $\dim $ $\mathbf{V}=n+m,n\geq
2,m\geq 1.$\ Local coordinates on $\mathbf{V}$ are labelled in the form $%
u=(x,y),$ or $u^{\alpha }=(x^{i},y^{a}),$ where indices $i,j,...=1,2,...,n$
are horizontal (h) ones and $a,b,...$ $=1,2,...,m$ are vertical (v) ones. We
follow our convention to use ''boldface'' symbols for nonholonomic spaces
and geometric objects on such spaces \cite{bejf,vnhrf1,vsgg,vrflg}. The
typical examples are those when $\mathbf{V}=TM$ is a tangent bundle, $%
\mathbf{V=E}$ is a vector bundle on $M,$ or $\mathbf{V}$ is a (semi--)
Riemann manifold, with prescribed local (nonintegrable) fibred structure.

In this work, a nonholonomic manifold $\mathbf{V}$ is considered to be
provided with a nonitegrable (nonhlonomic) distribution defining a nonlinear
connection (N--connection). This is equivalent to a Whitney sum of
conventional h-- and v--subspaces, $h\mathbf{V}$ and $v\mathbf{V},$%
\begin{equation}
T\mathbf{V}=h\mathbf{V\oplus }v\mathbf{V,}  \label{whit}
\end{equation}%
where $T\mathbf{V}$ is the tangent bundle. Such manifolds are called, in
brief, N--anholonomic (in literature, one uses two equivalent terms,
nonholonomic and anholonomic). Locally, a N--connection is defined by its
coefficients, $\mathbf{N=}\{N_{i}^{a}\},$ stated with respect to a local
coordinate basis, $\mathbf{N=}N_{i}^{a}(u)dx^{i}\otimes \partial /\partial
y^{a}.$ We can consider the class of linear connections when $%
N_{i}^{a}(u)=\Gamma _{ib}^{a}(x)y^{b}$ as a particular case.

N--connections are naturally considered in Finsler and Lagrange geometry %
\cite{ma1,ma2,bejf,vsgg}. They are related to (semi) spray configurations
\begin{equation}
\frac{dy^{a}}{d\varsigma }+2G^{a}(x,y)=0,  \label{ngeq}
\end{equation}%
of a curve $x^{i}(\varsigma )$ with parameter $0\leq \varsigma \leq
\varsigma _{0},$ when $y^{i}=dx^{i}/d\varsigma $ [spray configurations are
obtained for integrable equations]. For a regular Lagrangian $%
L(x,y)=L(x^{i},y^{a})$ modelled on $\mathbf{V,}$ when the Lagrange metric
(equivalently, Hessian)
\begin{equation}
\ ^{L}g_{ij}=\frac{1}{2}\frac{\partial ^{2}L}{\partial y^{i}\partial y^{j}}
\label{lm}
\end{equation}%
is not degenerate, i.e. $\det |g_{ij}|$ $\neq 0,$ one finds the fundamental
result (proof is a straightforward computation):

\begin{theorem}
\label{t1}For $4G^{j}=\ ^{L}g^{ij}\left( \frac{\partial ^{2}L}{\partial
y^{i}\partial x^{k}}y^{k}-\frac{\partial L}{\partial x^{i}}\right) ,$ with $%
\ ^{L}g^{ij}$ inverse to $\ ^{L}g_{ij},$ the ''nonlinear'' geodesic
equations (\ref{ngeq}) are equivalent to the Euler--Lagrange equations $%
\frac{d}{d\varsigma }\left( \frac{\partial L}{\partial y^{i}}\right) -\frac{%
\partial L}{\partial x^{i}}=0.$
\end{theorem}

Originally, the Lagrange geometry was elaborated on the tangent bundle $TM$
of a manifold $M,$ for a regular Lagrangian $L(x,y)$ following the methods
of Finsler geometry \cite{ma1,ma2} (Finsler configurations can be obtained
in a particular case when $L(x,y)=F^{2}(x,y)$ for a homogeneous fundamental
function $F(x,\lambda y)=\lambda F(x,y),\lambda \in \mathbb{R)}.$ Lagrange
and Finsler geometries can be also modelled on N--anholonomic manifolds \cite%
{bejf,vsgg} provided, for instance, with canonical N--connection structure
\begin{equation}
N_{i}^{a}=\frac{\partial G^{a}}{\partial y^{i}}.  \label{clnc}
\end{equation}

\begin{proposition}
A N--connection defines certain classes of nonholonomic preferred frames and
coframes,
\begin{eqnarray}
\mathbf{e}_{\alpha } &=&\left[ \mathbf{e}_{i}=\frac{\partial }{\partial x^{i}%
}-N_{i}^{a}(u)\frac{\partial }{\partial y^{a}},e_{b}=\frac{\partial }{%
\partial y^{b}}\right]  \label{dder} \\
\mathbf{e}^{\alpha } &=&[e^{i}=dx^{i},\mathbf{e}%
^{a}=dy^{a}+N_{i}^{a}(x,y)dx^{i}].  \label{ddif}
\end{eqnarray}
\end{proposition}

\begin{proof}
One computes the nontrivial nonholonomy coefficients $W_{ib}^{a}=\partial
N_{i}^{a}/\partial y^{b}$ and $W_{ij}^{a}=\Omega _{ji}^{a}=\mathbf{e}%
_{i}N_{j}^{a}-\mathbf{e}_{j}N_{i}^{a}$ (where $\Omega _{ji}^{a}$ are the
coefficients of the N--connection curvature) for
\begin{equation}
\left[ \mathbf{e}_{\alpha },\mathbf{e}_{\beta }\right] =\mathbf{e}_{\alpha }%
\mathbf{e}_{\beta }-\mathbf{e}_{\beta }\mathbf{e}_{\alpha }=W_{\alpha \beta
}^{\gamma }\mathbf{e}_{\gamma }.  \label{anhr}
\end{equation}%
$\square $
\end{proof}

One holds:

\begin{claim}
Any regular Lagrange mechanics $L(x,y)=L(x^{i},y^{a})$ modelled on $\mathbf{%
V,}\dim \mathbf{V}=2n,$ defines a canonical metric structure%
\begin{equation}
\ ^{L}\mathbf{g}=\ ^{L}g_{ij}(x,y)\left[ e^{i}\otimes e^{j}+\mathbf{e}%
^{i}\otimes \mathbf{e}^{j}\right] .  \label{m1}
\end{equation}
\end{claim}

\begin{proof}
For $\mathbf{V}=TM,$ the metric (\ref{m1}) is just the Sasaki lift of (\ref%
{lm}) on total space \cite{ma1,ma2}. In abstract form, such canonical
constructions can be performed similarly for any N--anholonomic manifold $%
\mathbf{V.}$ This approach to geometric mechanics follows from the fact that
the (semi) spray configurations are related to the N--connection structure
and defined both by the Lagrangian fundamental function and the
Euler--Lagrange equations, see Theorem \ref{t1} .$\square $
\end{proof}

\begin{definition}
A distinguished connection (d--connection) $\mathbf{D}$ on $\mathbf{V}$ is a
linear connection preserving under parallel transports the Whitney sum (\ref%
{whit}).
\end{definition}

In order to perform computations with d--connections we can use N--adapted
differential forms like $\mathbf{\Gamma }_{\ \beta }^{\alpha }=\mathbf{%
\Gamma }_{\ \beta \gamma }^{\alpha }\mathbf{e}^{\gamma }$ with the
coefficients defined with respect to (\ref{ddif}) and (\ref{dder}) and
parametrized $\mathbf{\Gamma }_{\ \alpha \beta }^{\gamma }=\left(
L_{jk}^{i},L_{bk}^{a},C_{jc}^{i},C_{bc}^{a}\right) .$ The torsion of a
d--connection is computed
\begin{equation}
\mathcal{T}^{\alpha }\doteqdot \mathbf{De}^{\alpha }=d\mathbf{e}^{\alpha
}+\Gamma _{\ \beta }^{\alpha }\wedge \mathbf{e}^{\beta }.  \label{tors}
\end{equation}%
Locally, it is characterized by (N--adapted) d--torsion coefficients
\begin{eqnarray}
T_{\ jk}^{i} &=&L_{\ jk}^{i}-L_{\ kj}^{i},\ T_{\ ja}^{i}=-T_{\ aj}^{i}=C_{\
ja}^{i},\ T_{\ ji}^{a}=\Omega _{\ ji}^{a},\   \notag \\
T_{\ bi}^{a} &=&-T_{\ ib}^{a}=\frac{\partial N_{i}^{a}}{\partial y^{b}}-L_{\
bi}^{a},\ T_{\ bc}^{a}=C_{\ bc}^{a}-C_{\ cb}^{a}.  \label{dtors}
\end{eqnarray}

\begin{theorem}
There is a unique canonical d--connection $\widehat{\mathbf{D}}=\{\widehat{%
\mathbf{\Gamma }}_{\ \alpha \beta }^{\gamma }=\left( \widehat{L}_{jk}^{i},%
\widehat{L}_{bk}^{a},\widehat{C}_{jc}^{i},\widehat{C}_{bc}^{a}\right) \}$
which is metric compatible with the Lagrange canonical metric structure, $%
\widehat{\mathbf{D}}\ \left( ^{L}\mathbf{g}\right) =0,$ and satisfies the
conditions $\widehat{T}_{\ jk}^{i}=$ $\widehat{T}_{\ bc}^{a}=0.$
\end{theorem}

\begin{proof}
It follows from explicit formulas (\ref{m1}) and (\ref{clnc}) and (\ref%
{candcon}). $\square $
\end{proof}

A geometric model of Lagrange mechanics can be elaborated in terms of
Riemannian geometry on $\mathbf{V,}$ as a noholonomic Riemann space, if we
chose the Levi--Civita connection $\nabla =\{\ _{\shortmid }\Gamma _{\
\alpha \beta }^{\gamma }\}$ defined uniquely by the Lagrange metric $\ ^{L}%
\mathbf{g}$ but such constructions are not adapted to the N--connection
splitting (\ref{whit}) induced by the (semi) spray Lagrange configuration.
In an equivalent form, such constructions can be adapted to the
N--connection structure if the canonical distinguished connection $\widehat{%
\mathbf{D}}=\{\widehat{\mathbf{\Gamma }}_{\ \alpha \beta }^{\gamma }\}$ is
considered. In this case, the geometric space is of Riemann--Cartan type,
with nontrivial torsion induced by the N--connection coefficients under
nonholonomic deformations of the frame structure.

\begin{conclusion}
Any regular Lagrange mechanics (Finsler ge\-ometry) can be modelled in two
equivalent canonical forms as a nonholonomic Riemann space or as a
N--anholonomic Riemann--Cartan space with the fundamental geometric objects
(metric and connection structures) defined by the fundamental Lagrange
(Finsler) function.
\end{conclusion}

Inverse statements when (semi) Riemannian metrics are modelled by certain
effective Lagrange structures and corresponding Ricci flows also hold true
but in such cases one has to work with models of generalized Lagrange
geometry, see \cite{vnhrf1,ma1,ma2}.

\begin{remark}
One considers different types of d--connec\-ti\-on structures in Finsler
geometry. For instance, there is an approach based on the so--called Chern
connection \cite{chern} which is not metric compatible and considered less
suitable for applications to standard models in modern physics, see
discussion in \cite{vsgg,vnhrf1}.
\end{remark}

For convenience, in Appendix, we outline the main formulas for the
connections $\nabla $ and $\widehat{\mathbf{D}}$ and their torsions,
curvature and Ricci tensors.

\section{The Perelman's Functionals on Lagrange and Finsler Spaces}

The Ricci flow equation was originally introduced by R. Hamilton \cite{ham1}
as an evolution equation
\begin{equation}
\frac{\partial g_{\alpha \beta }(\chi )}{\partial \chi }=-2\ _{\shortmid
}R_{\alpha \beta }(\chi )  \label{heq1}
\end{equation}%
for a set of Riemannian metrics $g_{\alpha \beta }(\chi )$ and corresponding
Ricci tensors $\ _{\shortmid }R_{\alpha \beta }(\chi )$ parametrized by a
real parameter $\chi .$\footnote{%
for our further purposes, on generalized Riemann--Finsler spaces, it is
convenient to use a different system of denotations than those considered by
R. Hamilton or Grisha Perelman on holonomic Riemannian spaces} The Ricci
flow theory is a branch of mathematics elaborated in connection to rigorous
study of topological and geometric properties of such equations and possible
applications in modern physics.

In the previous section, see also related details in our works \cite%
{vnhrf1,vsgg}, we proved that the Lagrange--Finsler geometries can be
modelled as constrained structures on N--anhonomic Riemannian spaces. We
concluded there that the Ricci flows of regular Lagrange systems (Finsler
metrics) can be described by usual Riemann gradient flows but subjected to
certain classes of nonholonomic constraints. It should be also noted that,
inversely, it is possible to extract from respective nonholonomic Riemannian
configurations the Lagrange or Finsler ones. Working with the canonical
d--connection $\widehat{\mathbf{D}},$ we get a Ricci tensor (\ref{dricci})
which, in general, is not symmetric but the metric (\ref{m1}) is symmetric.
In such cases, we are not able to derive the equation (\ref{heq1}) on
nonholonomic spaces in a self--consistent heuristic form following the
analogy of R. Hamilton's equations and the Einstein's equations. We
emphasize that one can be considered flows of nonholonomic Einstein spaces,
when $\widehat{\mathbf{R}}_{\alpha \beta }$ is symmetric (we investigated
such solutions in Refs. \cite{vijmpa1,vv1,vv2}), but more general classes of
solution of the Ricci equations with nonholonomic constraints would result
in nonsymmetric metrics, see discussions in Ref. \cite{vnhrf1}.

The Grisha Perelman's fundamental idea was to prove that the Ricci flow is
not only a gradient flow but also can be defined \ as a dynamical system on
the spaces of Riemannian metrics by introducing two Lyapunov type
functionals. In this section, we show how the constructions can be
generalized for N--anholonomic manifolds if we chose the connection $%
\widehat{\mathbf{D}}.$

The Perelman's functionals were introduced for Ricci flows of Riemannian
metrics. For the Levi--Civita connection defined by the Lagrange metric (\ref%
{m1}), are written in the form
\begin{eqnarray}
\ _{\shortmid }\mathcal{F}(L,f) &=&\int\limits_{\mathbf{V}}\left( \
_{\shortmid }R+\left| \nabla f\right| ^{2}\right) e^{-f}\ dV,  \label{pfrs}
\\
\ _{\shortmid }\mathcal{W}(L,f,\tau ) &=&\int\limits_{\mathbf{V}}\left[ \tau
\left( \ _{\shortmid }R+\left| \nabla f\right| \right) ^{2}+f-2n\right] \mu
\ dV,  \notag
\end{eqnarray}%
where $dV$ is the volume form of $\ ^{L}\mathbf{g,}$ integration is taken
over compact $\mathbf{V,}$  function $f$ is introduced in order to have the
possibility to consider gradient flows with different measures, see details
in \cite{gper1}, and $\ _{\shortmid }R$ is the scalar curvature computed for
$\nabla .$ For $\tau >0,$ we have $\int\nolimits_{\mathbf{V}}\mu dV=1$ when $%
\mu =\left( 4\pi \tau \right) ^{-n}e^{-f}.$%
\footnote{ In our works \cite{vrflg,vsgg}, we use left "up" and "low" indices as labels for some
geometric/ physical objects, for instance, in order  to emphasize that such values
are  induced by a  Lagrangian, ot defined by the Levi-Civita connection. }
The functional approach can be redefined for N--anholonomic manifolds:

\begin{claim}
For Lagrange spaces, the Perelman's functionals for the canonical
d--connection $\widehat{\mathbf{D}}$ are defined%
\begin{eqnarray}
\widehat{\mathcal{F}}(L,\widehat{f}) &=&\int\limits_{\mathbf{V}}\left(
R+S+\left| \widehat{\mathbf{D}}\widehat{f}\right| ^{2}\right) e^{-\widehat{f}%
}\ dV,  \label{npf1} \\
\widehat{\mathcal{W}}(L,\widehat{f},\tau ) &=&\int\limits_{\mathbf{V}}\left[
\widehat{\tau }\left( R+S+\left| ^{h}D\widehat{f}\right| +\left| ^{v}D%
\widehat{f}\right| \right) ^{2}+\widehat{f}-2n\right] \widehat{\mu }\ dV,
\label{npf2}
\end{eqnarray}%
where $dV$ is the volume form of $\ ^{L}\mathbf{g,}$ $R$ and $S$ are
respectively the h- and v--components of the curvature scalar of $\ \widehat{%
\mathbf{D}},$ see (\ref{sdccurv}), for $\ \widehat{\mathbf{D}}_{\alpha
}=(D_{i},D_{a}),$ or $\widehat{\mathbf{D}}=(\ ^{h}D,\ ^{v}D),$ $\left|
\widehat{\mathbf{D}}\widehat{f}\right| ^{2}=\left| ^{h}D\widehat{f}\right|
^{2}+\left| ^{v}D\widehat{f}\right| ^{2},$ and $\widehat{f}$ satisfies $%
\int\nolimits_{\mathbf{V}}\widehat{\mu }dV=1$ for $\widehat{\mu }=\left(
4\pi \tau \right) ^{-n}e^{-\widehat{f}}$ and $\tau >0.$
\end{claim}

\begin{proof}
We can redefine equivalently the formulas (\ref{pfrs}) for some $\widehat{f}$
and $f$ (which can be a non--explicit relation) when
\begin{equation*}
\left( \ _{\shortmid }R+\left| \nabla f\right| ^{2}\right) e^{-f}=\left(
R+S+\left| ^{h}D\widehat{f}\right| ^{2}+\left| ^{v}D\widehat{f}\right|
^{2}\right) e^{-\widehat{f}}\ +\Phi
\end{equation*}%
and re--scale the parameter $\tau \rightarrow \widehat{\tau }$ to have
\begin{equation*}
[\tau (\ _{\shortmid }R+\left| \nabla f\right| )^{2}+f-2n]\mu =[\widehat{%
\tau }(R+S+\left| ^{h}D\widehat{f}\right| +\left| ^{v}D\widehat{f}\right|
)^{2}+\widehat{f}-2n]\widehat{\mu }+\Phi _{1}
\end{equation*}%
for some $\Phi $ and $\Phi _{1}$ for which $\int\limits_{\mathbf{V}}\Phi
dV=0 $ and $\int\limits_{\mathbf{V}}\Phi _{1}dV=0.$ $\square $
\end{proof}

Elaborating a N--adapted variational calculus, we shall consider both
variations in the so--called h-- and v--subspaces as defined by the
decompositions (\ref{whit}). We write, for simplicity, $g_{ij}=\ ^{L}g_{ij}$
and consider the h--variation $^{h}\delta g_{ij}=v_{ij},$ the v--variation $%
^{v}\delta g_{ab}=v_{ab},$ for a fixed N--connection structure in (\ref{m1}%
), and $^{h}\delta \widehat{f}=\ ^{h}f,$ $^{v}\delta \widehat{f}=\ ^{v}f$

\begin{lemma}
\label{lem1}The first N--adapted variations of (\ref{npf1}) are given by
\begin{eqnarray}
&&\delta \widehat{\mathcal{F}}(v_{ij},v_{ab},\ ^{h}f,\ ^{v}f)=  \label{vnpf1}
\\
&&\int\limits_{\mathbf{V}}\{[-v_{ij}(R_{ij}+D_{i}D_{j}\widehat{f})+(\frac{\
^{h}v}{2}-\ ^{h}f)\left( 2\ ^{h}\Delta \widehat{f}-|\ ^{h}D\ \widehat{f}%
|\right) +R]  \notag \\
&&+[-v_{ab}(R_{ab}+D_{a}D_{b}\widehat{f})+(\frac{\ ^{v}v}{2}-\ ^{v}f)\left(
2\ ^{v}\Delta \widehat{f}-|\ ^{v}D\ \widehat{f}|\right) +S]\}e^{-\widehat{f}%
}dV  \notag
\end{eqnarray}%
where $\widehat{\Delta }=$ $\ ^{h}\Delta +\ ^{v}\Delta ,$ $^{h}\Delta
=D_{i}D^{i},$ $^{v}\Delta =D_{a}D^{a},$ $\ ^{h}v=g^{ij}v_{ij},\
^{v}v=g^{ab}v_{ab}.$
\end{lemma}

\begin{proof}
It is a N--adapted calculus similar to that for Perelman's Lemma in \cite%
{gper1}. We omit details given, for instance, in the proof from \cite{caozhu}%
, see there Lemma 1.5.2, but we note that if such computations are performed
on a N--anholonomic manifold, the canonical d--connection results in
formulas (\ref{dricci}), for the Ricci curvature and (\ref{sdccurv}), for
the scalar curvature of $\widehat{\mathbf{D}}.$ It should be emphasized that
because we consider that variations of a symmetric metric, $^{h}\delta
g_{ij}=v_{ij}$ and $^{v}\delta g_{ab}=v_{ab},$ are considered independently
on h-- and v--subspaces and supposed to be also symmetric, we get in (\ref%
{vnpf1}) only the symmetric coefficients $R_{ij}$ and $R_{ab}$ but not $%
R_{ai}$ and $R_{ia}.$ Admitting nonsymmetric variations of metrics, we would
obtain certain terms in $\delta \widehat{\mathcal{F}}(v_{ij},v_{ab},\
^{h}f,\ ^{v}f)$ defined by the nonsymmetric components of the Ricci tensor
for $\widehat{\mathbf{D}}.$ In this work, we try to keep our constructions
on Riemannian spaces, even they are provided with N--anholonomic
distributions, and avoid to consider the so--called Lagrange--Eisenhart, or
Finsler--Eisenhart, geometry analyzed, for instance in Chapter 8 of
monograph \cite{ma1} (for nonholonomic Ricci flows, we discuss the problem
in \cite{vnhrf1}). $\square $
\end{proof}

\section{Evolution Equations for Lagrange Systems}

The normalized (holonomic) Ricci flows, see details in Refs. \cite%
{gper1,caozhu,kleiner,rbook}, with respect to the coordinate base $\partial
_{\underline{\alpha }}=\partial /\partial u^{\underline{\alpha }},$ are
described by the equations
\begin{equation}
\frac{\partial }{\partial \chi }g_{\underline{\alpha }\underline{\beta }%
}=-2\ _{\shortmid }R_{\underline{\alpha }\underline{\beta }}+\frac{2r}{5}g_{%
\underline{\alpha }\underline{\beta }},  \label{feq}
\end{equation}%
where the normalizing factor $r=\int \ _{\shortmid }RdV/dV$ is introduced in
order to preserve the volume $V.$ \footnote{%
we underlined the indices with respect to the coordinate bases in order to
distinguish them from those defined with respect to the 'N--elongated' local
bases (\ref{dder}) and (\ref{ddif})} We note that here we use the Ricci
tensor $\ _{\shortmid }R_{\underline{\alpha }\underline{\beta }}$ and scalar
curvature $\ _{\shortmid }R=g^{\underline{\alpha }\underline{\beta }}\
_{\shortmid }R_{\underline{\alpha }\underline{\beta }}$ computed for the
connection $\nabla .$ The coefficients $g_{\underline{\alpha }\underline{%
\beta }}$ are those for for a family of metrics (\ref{m1}), $\ ^{L}\mathbf{g}%
\left( \chi \right) $ ,\ rewritten with respect to the coordinate basis, $\
^{L}\mathbf{g}\left( \chi \right) =g_{\underline{\alpha }\underline{\beta }%
}(\chi )du^{\underline{\alpha }}\otimes du^{\underline{\beta }},$ where%
\begin{equation}
\ g_{\underline{\alpha }\underline{\beta }}(\chi )=\left[
\begin{array}{cc}
\underline{g}_{ij}(\chi )=\ g_{ij}+~N_{i}^{a}~N_{j}^{b}\ g_{ab} & ~%
\underline{g}_{ib}(\chi )=~N_{i}^{e}\ g_{be} \\
~\underline{g}_{aj}(\chi )=N_{i}^{e}\ g_{be} & \ \underline{g}_{ab}(\chi )%
\end{array}%
\right] ,  \label{qel}
\end{equation}%
for $g_{\underline{\alpha }\underline{\beta }}(\chi )=\underline{g}_{\alpha
\beta }(\chi ),$ when $g_{ij}(\chi )=\ ^{L}g_{ij}(\chi ,u),g_{ab}(\chi )=\
^{L}g_{ab}(\chi ,u)$ and $~N_{i}^{a}(\chi )~=N_{i}^{a}(\chi ,u)$ defined
from a set of Lagrangians $L(\chi ,u),$ respectively by formulas (\ref{lm})
and (\ref{clnc}).

With respect to the N--adapted frames (\ref{dder}) and (\ref{ddif}), when
\begin{equation*}
\mathbf{e}_{\alpha }(\chi )=\mathbf{e}_{\alpha }^{\ \underline{\alpha }%
}(\chi )\ \partial _{\underline{\alpha }}\mbox{\ and \ }\mathbf{e}^{\alpha
}(\chi )=\mathbf{e}_{\ \underline{\alpha }}^{\alpha }(\chi )du^{\underline{%
\alpha }},
\end{equation*}%
the frame transforms are respectively parametrized in the form%
\begin{eqnarray}
\mathbf{e}_{\alpha }^{\ \underline{\alpha }}(\chi ) &=&\left[
\begin{array}{cc}
e_{i}^{\ \underline{i}}=\delta _{i}^{\underline{i}} & e_{i}^{\ \underline{a}%
}=N_{i}^{b}(\chi )\ \delta _{b}^{\underline{a}} \\
e_{a}^{\ \underline{i}}=0 & e_{a}^{\ \underline{a}}=\delta _{a}^{\underline{a%
}}%
\end{array}%
\right] ,  \label{ft} \\
\mathbf{e}_{\ \underline{\alpha }}^{\alpha }(\chi ) &=&\left[
\begin{array}{cc}
e_{\ \underline{i}}^{i}=\delta _{\underline{i}}^{i} & e_{\ \underline{i}%
}^{b}=-N_{k}^{b}(\chi )\ \delta _{\underline{i}}^{k} \\
e_{\ \underline{a}}^{i}=0 & e_{\ \underline{a}}^{a}=\delta _{\underline{a}%
}^{a}%
\end{array}%
\right] ,  \notag
\end{eqnarray}%
where $\delta _{\underline{i}}^{i}$ is the Kronecher symbol, the Ricci flow
equations (\ref{feq}) are
\begin{eqnarray}
&&\frac{\partial }{\partial \chi }g_{ij}=2[N_{i}^{a}N_{j}^{b}(\ _{\shortmid
}R_{ab}-\lambda g_{ab})-\ _{\shortmid }R_{ij}+\lambda g_{ij}]-g_{cd}\frac{%
\partial }{\partial \chi }(N_{i}^{c}N_{j}^{d}),  \label{eq1} \\
&&\frac{\partial }{\partial \chi }g_{ab}=-2\ _{\shortmid }R_{ab}+2\lambda
g_{ab},\   \label{eq2} \\
&&\frac{\partial }{\partial \chi }(N_{j}^{e}\ g_{ae})=-2\ _{\shortmid
}R_{ia}+2\lambda N_{j}^{e}\ g_{ae},  \label{eq3}
\end{eqnarray}%
where $\lambda =r/5$ and the metric coefficients are defined by the ansatz (%
\ref{qel}).

If $\nabla \rightarrow \widehat{\mathbf{D}},$ we have to change $_{\shortmid
}R_{\alpha \beta }\rightarrow \widehat{\mathbf{R}}_{\alpha \beta }$ in (\ref%
{eq1})--(\ref{eq3}). The N--adapted evolution equations for Ricci flows of
symmetric metrics, with respect to local coordinate frames, are written
\begin{eqnarray}
\frac{\partial }{\partial \chi }g_{ij} &=&2[N_{i}^{a}N_{j}^{b}(\underline{%
\widehat{R}}_{ab}-\lambda g_{ab})-\underline{\widehat{R}}_{ij}+\lambda
g_{ij}]-g_{cd}\frac{\partial }{\partial \chi }(N_{i}^{c}N_{j}^{d}),
\label{e1} \\
\frac{\partial }{\partial \chi }g_{ab} &=&-2(\underline{\widehat{R}}%
_{ab}-\lambda g_{ab}),\   \label{e2} \\
\ \widehat{R}_{ia} &=&0\mbox{ and }\ \widehat{R}_{ai}=0,  \label{e3}
\end{eqnarray}%
where the Ricci coefficients $\underline{\widehat{R}}_{ij}$ and $\underline{%
\widehat{R}}_{ab}$ are computed with respect to coordinate coframes, being
frame transforms (\ref{ft}) of the corresponding formulas (\ref{dricci})
defined with respect to N--adapted frames. The equations (\ref{e3})
constrain the nonholonomic Ricci flows to result in symmetric metrics.

The aim of this section is to prove that equations of type (\ref{e1}) and (%
\ref{e2}) can be derived from the Perelman's N--adapted functionals (\ref%
{npf1}) and (\ref{npf2}) (for simplicity, we shall not consider the
normalized term and put $\lambda =0).$

\begin{definition}
A metric $\ ^{L}\mathbf{g}$ generated by a regular Lagrangian $L$ evolving
by the (nonholonomic) Ricci flow is called a (nonholonomic) breather if for
some $\chi _{1}<\chi _{2}$ and $\alpha >0$ the metrics $\alpha \ ^{L}\mathbf{%
g(}\chi _{1}\mathbf{)}$ and $\alpha \ ^{L}\mathbf{g(}\chi _{2}\mathbf{)}$
differ only by a diffeomorphism (in the N--anholonomic case, preserving the
Whitney sum (\ref{whit})). The cases $\alpha =,<,>1$ define correspondingly
the steady, shrinking and expanding breathers (for N--anholonomic manifolds,
one can be the situation when, for instance, the h--component of metric is
steady but the v--component is shrinking).
\end{definition}

Clearly, the breather properties depend on the type of connections are used
for definition of Ricci flows.

Following a N--adapted variational calculus for $\widehat{\mathcal{F}}(L,%
\widehat{f}),$ see Lemma \ref{lem1}, with Laplacian $\widehat{\Delta }$ and
h- and v--components of the Ricci tensor, $\widehat{R}_{ij}$ and $\widehat{S}%
_{ij},$ defined by $\widehat{\mathbf{D}}$ and considering parameter $\tau
(\chi ),$ $\partial \tau /\partial \chi =-1,$ we prove

\begin{theorem}
The Ricci flows of regular Lagrange mechanical systems are characterized by
evolution equations
\begin{equation*}
\frac{\partial \underline{g}_{ij}}{\partial \chi }=-2\underline{\widehat{R}}%
_{ij},\ \frac{\partial \underline{g}_{ab}}{\partial \chi }=-2\underline{%
\widehat{R}}_{ab},\ \frac{\partial \widehat{f}}{\partial \chi }=-\widehat{%
\Delta }\widehat{f}+\left| \widehat{\mathbf{D}}\widehat{f}\right| ^{2}-R-S
\end{equation*}%
and the property that, for constant $\int\limits_{\mathbf{V}}e^{-\widehat{f}%
}dV,$
\begin{equation*}
\frac{\partial }{\partial \chi }\widehat{\mathcal{F}}(\ ^{L}\mathbf{g(\chi ),%
}\widehat{f}(\chi ))=2\int\limits_{\mathbf{V}}\left[ |\widehat{R}%
_{ij}+D_{i}D_{j}\widehat{f}|^{2}+|\widehat{R}_{ab}+D_{a}D_{b}\widehat{f}|^{2}%
\right] e^{-\widehat{f}}dV.
\end{equation*}
\end{theorem}

\begin{proof}
For Riemannian spaces, a proof was proposed by G. Perelman \cite{gper1}
(details of the proof are given for the connection $\nabla $ in the
Proposition 1.5.3 of \cite{caozhu}, they can be similarly reproduced for the
canonical d--connection $\widehat{\mathbf{D}}$). For N--anholonomic spaces,
we changed the status of such statements to a Theorem because for
nonholonomic configurations there are not alternative ways of definition
Ricci flow equations in N--adapted form following two different, heuristic
and functional, approaches. The functional variant became the unique
possibility for a rigorous proof containing N--adapted calculus. Finally, we
note that for the Levi--Civita connection the functional $\mathcal{F}$ is
nondecreasing in time and the monotonicity is strict unless we are on a
steady gradient soliton (see, for instance, Ref. \cite{caozhu} for details
on solitonic solutions and Ricci flows). This property sure depend on the
type of connection is used and how solitons are defined. We shall not use it
in this works and omit such considerations.$\square $
\end{proof}

The priority of the N--adapted calculus for the canonical d--connection $%
\widehat{\mathbf{D}}=(^{h}D,\ ^{v}D)$ is that from formal point of view we
work as in the case with the connection $\nabla $ but have to dub the
results for the h-- and v--components and redefine them with respect
N--adapted bases. This analogy holds true for all (generalized) and Lagrange
and Finsler metrics because $\widehat{\mathbf{D}}$ is metric compatible and
uniquely defined by the coefficients of $\ ^{L}\mathbf{g,}$ similarly to $%
\nabla .$ It should be noted that even a closed formal analogy of formulas
exist, the evolution equations, their solutions, and related geometrical and
fundamental objects are different because $\widehat{\mathbf{D}}\neq \nabla .$
Following this property, we can formulate (the reader may check that its
statements and proofs consist a N--adapted modification of Proposition 1.5.8
in \cite{caozhu} containing the details of the original result from \cite%
{gper1}):

\begin{theorem}
\label{theveq}If a regular Lagrange (Finsler) metric $\ ^{L}\mathbf{g}(\chi
) $ and functions $\widehat{f}(\chi )$ and $\widehat{\tau }(\chi )$ evolve
for $\frac{\partial \widehat{\tau }}{\partial \chi }=-1$ and constant $%
\int\limits_{\mathbf{V}}(4\pi \widehat{\tau })^{-n}e^{-\widehat{f}}dV,$ as
solutions of the system
\begin{equation*}
\frac{\partial \underline{g}_{ij}}{\partial \chi }=-2\underline{\widehat{R}}%
_{ij},\ \frac{\partial \underline{g}_{ab}}{\partial \chi }=-2\underline{%
\widehat{R}}_{ab},\ \frac{\partial \widehat{f}}{\partial \chi }=-\widehat{%
\Delta }\widehat{f}+\left| \widehat{\mathbf{D}}\widehat{f}\right| ^{2}-R-S+%
\frac{n}{\widehat{\tau }},
\end{equation*}
one holds the equality
\begin{eqnarray*}
\frac{\partial }{\partial \chi }\widehat{\mathcal{W}}(\ ^{L}\mathbf{g}(\chi )%
\mathbf{,}\widehat{f}(\chi ),\widehat{\tau }(\chi )) &=&2\int\limits_{%
\mathbf{V}}\widehat{\tau }[|\widehat{R}_{ij}+D_{i}D_{j}\widehat{f}-\frac{1}{2%
\widehat{\tau }}g_{ij}|^{2}+ \\
&&|\widehat{R}_{ab}+D_{a}D_{b}\widehat{f}-\frac{1}{2\widehat{\tau }}%
g_{ab}|^{2}](4\pi \widehat{\tau })^{-n}e^{-\widehat{f}}dV.
\end{eqnarray*}
\end{theorem}

For the Levi--Civita connection $\nabla ,$ the functional $\ _{\shortmid }%
\mathcal{W}(\ ^{L}\mathbf{g}(\chi )\mathbf{,}f(\chi ),\tau (\chi ))$ is
nondecreasing in time and the monotonicity is strict unless we are on a
shrinking gradient soliton. Similar properties can be formulated in
N--adapted form, but it is not obvious if some of them hold true for $\nabla
$ they will be preserved for $\nabla \rightarrow \widehat{\mathbf{D}}.$

The Lagrange--Ricci flows are are characterized by the evolutions of
preferred N--adapted frames (\ref{ft}) (see proof in \cite{vnhrf1}):

\begin{corollary}
The evolution, for all time $\tau \in \lbrack 0,\tau _{0}),$ of preferred
frames on a Lagrange space $\ \mathbf{e}_{\alpha }(\tau )=\ \mathbf{e}%
_{\alpha }^{\ \underline{\alpha }}(\tau ,u)\partial _{\underline{\alpha }}$
is defined by the coefficients
\begin{equation*}
\ \mathbf{e}_{\alpha }^{\ \underline{\alpha }}(\tau ,u)=\left[
\begin{array}{cc}
\ e_{i}^{\ \underline{i}}(\tau ,u) & ~N_{i}^{b}(\tau ,u)\ e_{b}^{\
\underline{a}}(\tau ,u) \\
0 & \ e_{a}^{\ \underline{a}}(\tau ,u)%
\end{array}%
\right] ,\
\end{equation*}%
with $\ \ ^{L}g_{ij}(\tau )=\ e_{i}^{\ \underline{i}}(\tau ,u)\ e_{j}^{\
\underline{j}}(\tau ,u)\eta _{\underline{i}\underline{j}},$ where $\eta _{%
\underline{i}\underline{j}}=diag[\pm 1,...\pm 1]$ states the signature of $\
^{L}g_{\underline{\alpha }\underline{\beta }}^{[0]}(u),$ is given by
equations
\begin{eqnarray*}
\frac{\partial }{\partial \tau }\ e_{\alpha }^{\ \underline{\alpha }} &=&\
^{L}g^{\underline{\alpha }\underline{\beta }}~_{\shortmid }R_{\underline{%
\beta }\underline{\gamma }}~\ e_{\alpha }^{\ \underline{\gamma }},\ \ %
\mbox{\ for the Levi-Civita connection }; \\
\frac{\partial }{\partial \tau }\ e_{\alpha }^{\ \underline{\alpha }} &=&\
^{L}g^{\underline{\alpha }\underline{\beta }}~\widehat{R}_{\underline{\beta }%
\underline{\gamma }}~\ e_{\alpha }^{\ \underline{\gamma }},\ \
\mbox{\ for
the canonical d--connection }.
\end{eqnarray*}
\end{corollary}

It should be emphasized that it would be a problem to prove directly the results of this section for Ricci flows  of Finsler spaces with metric noncompatible d--connections like in Ref. \cite{chern}. Nevertheless, our proofs can be generalized also for nonmetric Lagrange--Finsler configurations if the nonmetricity is completely defined by the coefficients of the d--metric and N--connection structures. In such cases, we can prove the theorems and consequences as for metric compatible cases (for the Levi--Civita connection and/or Cartan d--connection) and then to distort the formulas in unique forms using corresponding deformation tensors.

\section{Statistical Analogy for Lagrange--Finsler \newline
Spaces and Ricci Flows}

Grisha Perelman showed that the functional $\ _{\shortmid }\mathcal{W}$ is
in a\ sense analogous to minus entropy \cite{gper1}. We show that this
property holds true for nonholonomic Ricci flows which provides a
statistical model for regular Lagrange (Finsler) systems.

The partition function $Z=\int \exp (-\beta E)d\omega (E)$ for the canonical
ensemble at temperature $\beta ^{-1}$ is defined by the measure taken to be
the density of states $\omega (E).$ The thermodynamical values are computed
in the form: the average energy, $<E>=-\partial \log Z/\partial \beta ,$ the
entropy $S=\beta <E>+\log Z$ and the fluctuation $\sigma =<\left(
E-<E>\right) ^{2}>=\partial ^{2}\log Z/\partial \beta ^{2}.$

Let us suppose that a set of regular mechanical systems with Lagrangians $L(%
\widehat{\tau },x,y)$ is described by respective metrics $^{L}\mathbf{g}(%
\widehat{\tau })$ and N--connection $N_{i}^{a}(\widehat{\tau })$ and related
canonical linear connections $\nabla (\widehat{\tau })$ and $\widehat{%
\mathbf{D}}(\widehat{\tau })$ subjected to the conditions of Theorem \ref%
{theveq}. One holds

\begin{theorem}
Any family of regular Lagrange (Finsler) geometries satisfying the evolution
equations for the canonical d--connection is characterized by thermodynamic
values
\begin{eqnarray*}
&<&\widehat{E}>\ =-\widehat{\tau }^{2}\int\limits_{\mathbf{V}}\left(
R+S+\left| ^{h}D\widehat{f}\right| ^{2}+\left| ^{v}D\widehat{f}\right| ^{2}-%
\frac{n}{\widehat{\tau }}\right) \widehat{\mu }\ dV, \\
\widehat{S} &=&-\int\limits_{\mathbf{V}}\left[ \widehat{\tau }\left(
R+S+\left| ^{h}D\widehat{f}\right| ^{2}+\left| ^{v}D\widehat{f}\right|
^{2}\right) +\widehat{f}-2n\right] \widehat{\mu }\ dV, \\
\widehat{\sigma } &=&2\ \widehat{\tau }^{4}\int\limits_{\mathbf{V}}\left[ |%
\widehat{R}_{ij}+D_{i}D_{j}\widehat{f}-\frac{1}{2\widehat{\tau }}%
g_{ij}|^{2}+|\widehat{R}_{ab}+D_{a}D_{b}\widehat{f}-\frac{1}{2\widehat{\tau }%
}g_{ab}|^{2}\right] \widehat{\mu }\ dV.
\end{eqnarray*}
\end{theorem}

\begin{proof}
It follows from a straightforward computation for $\widehat{Z}=$ $\exp
\{\int\nolimits_{\mathbf{V}}$ $[-\widehat{f}+n]\widehat{\mu }dV\}.$ We note
that similar values $<\ _{\shortmid }E>,\ _{\shortmid }S$ and $\ _{\shortmid
}\sigma $ can be computed for the Levi--Civita connection $\nabla $ also
defined for the metric $^{L}\mathbf{g,}$ see functionals (\ref{pfrs}). $%
\square $
\end{proof}

This results in

\begin{corollary}
A N--anholonomic Lagrange (Finsler) model defined by the canonical
d--connection $\widehat{\mathbf{D}}$ is thermodynamically more (less,
equivalent) convenient than a similar one defined by the Levi--Civita
connection $\nabla $ if $\ \widehat{S}<\ _{\shortmid }S$ ($\widehat{S}>\
_{\shortmid }S,\widehat{S}=\ _{\shortmid }S$).
\end{corollary}

Following this Corollary, we conclude that such models are positively
equivalent for integrable N--anholonomic structures with vanishing
distorsion tensor (see formulas (\ref{cdeft}) and (\ref{cdeftc})). For such
holonomic structures, the anholonomy coefficients $W_{\alpha \beta }^{\gamma
}$ (\ref{anhr}) are zero and we can work only with the Levi--Civita
connection. There are necessary explicit computations of the thermodynamical
values for different classes of exact solutions of nonholonomic Ricci flow
equations \cite{vijmpa1,vv1,vv2} or of the Einstein equations with
nonholonomic/ noncommutative variables \cite{vsgg} in order to conclude
which configurations are thermodynamically more convenient for
N--anholonomic or (pseudo) Riemannian configurations. In certain cases, some
constrained (Finsler like, or more general) configurations may be more
optimal than the Levi--Civita ones.

Finally, we would like to mention that there were elaborated alternative
approaches to geometric and non--equilibrium thermodynamics, locally
anisotropic kinetics and kinetic processes elaborated in terms of Riemannian
and Finsler like objects on phase and thermodynamic spaces, see reviews of
results and bibliography in Ref. \cite{rup,mrug,salb,aim,rap,vap} and
Chapter 6 from \cite{vsgg}. Those models are not tailor-made for Ricci flows of
geometric objects and seem not to be related to the statistical
thermodynamics of metrics and connections which can be derived from (an)
holonomic Perelman's functionals.  In a more general context, the ''Ricci
flow thermodynamics'' seem to be related to ''non-extensive" Tsallis
statistics which is valid for non--equilibrium cases and is considered to
be ''more fundamental'' than the equilibrium Boltzman--Gibbs statistics, see
Ref. \cite{tsallis} and references therein.

\section{Conclusion and Discussion}

In this paper, we have introduced an extension of Prelman's functional
approach to Ricci flows \cite{gper1} in order to derive in canonical form
the evolution equations for Lagrange and Finsler geometries and formulate a
statistical analogy of regular mechanical systems. This scheme is of
practical applicability to the problem of the definition of the most optimal
geometric and topological configurations in geometric mechanics and
analogous models of field interactions. In this context, we elaborate a new
direction to geometrization of Lagrange systems following the theory of
nonholonomic Ricci flows and generalized Riemann--Cartan and
Lagrange--Finsler spaces equipped with compatible metric, nonlinear
connection and linear connection structures \cite{vijmpa1,vv1,vv2,vsgg}.

Since the initial works on Ricci flows \cite{ham1,caozhu,kleiner,rbook}, the
problem of definition of evolution equations of fundamental geometric
objects was treated in a heuristic form following certain analogy with the
original 'proof' of the Einstein equations when a symmetric Ricci tensor was
set to be proportional to a 'simple' and physically grounded combination of
coordinate/parametric derivatives of metric coefficients. In our works \cite%
{vsgg}, we proved that Finsler like geometries can be modelled by preferred
nonholonomic frame structures even as exact solutions in the Einstein and
string gravity and has analogous interpretations in terms of geometric
objects on generalized Lagrange spaces and nonholonomic manifolds \cite%
{ma1,ma2,vsgg}. Then, it was shown that flows of Lagrange--Finsler
geometries can be extracted from flows of Riemannian metrics by imposing
certain classes of nonholonomic constraints and deformations of the frame
and linear connection structures \cite{vnhrf1}.

In order to derive the first results on Lagrange--Ricci, or Finsler--Ricci
flows, in a form more familiar to researches skilled in geometric analysis
and Riemannian geometry, we worked in the bulk with the Levi--Civita
connection for Lagrange, or Finsler, metrics and then sketched how the
results can be redefined in terms of the canonical connections for 'locally
isotropic' geometries. The advantage of the Perelman's approach to the Ricci
flow theory is that it can easily be reformulated for a covariant calculus
adapted to the nonlinear connection structure which is of crucial importance
in generalized Riemann--Finsler geometry. For such geometries, the
functional methods became a strong tool both for rigorous proofs of the
nonholonomic evolution equations and formulating new alternative statistical
models for regular Lagrange systems.

The two approaches are complementary in the following sense: the functional
scheme gives more rigorous results when the type of geometric structures are
prescribed and the holonomic or nonholonomic Ricci flows and the related
statistical/thermodynamical models are constructed in the same class of
geometries, whereas the heuristic ideas and formulas are best adapted for
flow transitions from one type of geometries to another ones (for instance,
from Finsler configurations to Riemannian ones, and inversely).

The next challenge in our program on nonholonomic Ricci flows and
applications is to formulate a functional formalism for general nonholonomic
manifolds in a form, when various type of nonholonomic Clifford, algebroid,
noncommutative, solitonic ... structures can be extracted from flows of
'Riemannian' geometrical objects by imposing the corresponding classes of
nonholonomic constraints and deformations of geometric objects. We discuss
such results and provide a more detailed list of references on Ricci flows
and applications to modern classical and quantum physics in our recent work %
\cite{vijmpa1,vv1,vv2,vrf2,vijtp,vncrf} (see references therein).

\vskip5pt

\textbf{Acknowledgements }

The author thanks Professors Mihai Anastasiei and Aurel Bejancu for valuable
support and important references on the geometry of Lagrange--Finsler spaces
and nonholonomic manifolds. He is also grateful to the first referee for an
important remark about ''Tsallis statistics'' and its possible connections
to Ricci flow thermodynamics.

\setcounter{equation}{0} \renewcommand{\theequation}
{A.\arabic{equation}} \setcounter{subsection}{0}
\renewcommand{\thesubsection}
{A.\arabic{subsection}}

\appendix

\section{Appendix}

One exists a minimal extension of the Levi--Civita connection $\nabla $ to a
canonical d--connection $\widehat{\mathbf{D}}$ which is defined only the
coefficients of Lagrange metric $^{L}\mathbf{g}$ (\ref{m1}) and canonical
nonlinear connection $N_{i}^{a}$ (\ref{clnc}) which is also metric
compatible, with $\widehat{T}_{\ jk}^{i}=0$ and $\widehat{T}_{\ bc}^{a}=0,$
but $\widehat{T}_{\ ja}^{i},\widehat{T}_{\ ji}^{a}$ and $\widehat{T}_{\
bi}^{a}$ are not zero, see (\ref{dtors}). The coefficient $\widehat{\mathbf{%
\Gamma }}_{\ \alpha \beta }^{\gamma }=\left( \widehat{L}_{jk}^{i},\widehat{L}%
_{bk}^{a},\widehat{C}_{jc}^{i},\widehat{C}_{bc}^{a}\right) $ of this
connection, with respect to the N--adapted frames, are computed:
\begin{eqnarray}
\widehat{L}_{jk}^{i} &=&\frac{1}{2}g^{ir}\left(
e_{k}g_{jr}+e_{j}g_{kr}-e_{r}g_{jk}\right) ,  \label{candcon} \\
\widehat{L}_{bk}^{a} &=&e_{b}(N_{k}^{a})+\frac{1}{2}g^{ac}\left(
e_{k}g_{bc}-g_{dc}\ e_{b}N_{k}^{d}-g_{db}\ e_{c}N_{k}^{d}\right) ,  \notag \\
\widehat{C}_{jc}^{i} &=&\frac{1}{2}g^{ik}e_{c}g_{jk},\ \widehat{C}_{bc}^{a}=%
\frac{1}{2}g^{ad}\left( e_{c}g_{bd}+e_{c}g_{cd}-e_{d}g_{bc}\right) ,  \notag
\end{eqnarray}%
where, for simplicity, we write $g_{jr}$ and $g_{bd}$ without label ''L'' we
used for Hessian $\ ^{L}g_{ij}$ (\ref{lm}).

The Levi--Civita linear connection $\bigtriangledown =\{\ _{\shortmid
}\Gamma _{\beta \gamma }^{\alpha }\},$ uniquely defined by the conditions $%
~\ _{\shortmid }\mathcal{T}=0$ and $\bigtriangledown g=0,$ is not adapted to
the distribution (\ref{whit}). Let us parametrize the coefficients in the
form
\begin{eqnarray*}
_{\shortmid }\Gamma _{\beta \gamma }^{\alpha } &=&\left( \ _{\shortmid
}L_{jk}^{i},\ _{\shortmid }L_{jk}^{a},\ _{\shortmid }L_{bk}^{i},\
_{\shortmid }L_{bk}^{a},\ _{\shortmid }C_{jb}^{i},\ _{\shortmid
}C_{jb}^{a},\ _{\shortmid }C_{bc}^{i},\ _{\shortmid }C_{bc}^{a}\right) , \\
\bigtriangledown _{\mathbf{e}_{k}}(\mathbf{e}_{j}) &=&\ _{\shortmid
}L_{jk}^{i}\mathbf{e}_{i}+\ _{\shortmid }L_{jk}^{a}e_{a},\ \bigtriangledown
_{\mathbf{e}_{k}}(e_{b})=\ _{\shortmid }L_{bk}^{i}\mathbf{e}_{i}+\
_{\shortmid }L_{bk}^{a}e_{a}, \\
\bigtriangledown _{e_{b}}(\mathbf{e}_{j}) &=&\ _{\shortmid }C_{jb}^{i}%
\mathbf{e}_{i}+\ _{\shortmid }C_{jb}^{a}e_{a},\ \bigtriangledown
_{e_{c}}(e_{b})=\ _{\shortmid }C_{bc}^{i}\mathbf{e}_{i}+\ _{\shortmid
}C_{bc}^{a}e_{a}.
\end{eqnarray*}%
It is convenient to express
\begin{equation}
\ _{\shortmid }\Gamma _{\ \alpha \beta }^{\gamma }=\widehat{\mathbf{\Gamma }}%
_{\ \alpha \beta }^{\gamma }+\ _{\shortmid }Z_{\ \alpha \beta }^{\gamma }
\label{cdeft}
\end{equation}%
where the explicit components of distorsion tensor $\ _{\shortmid }Z_{\
\alpha \beta }^{\gamma }$ \ are computed
\begin{eqnarray}
\ _{\shortmid }Z_{bk}^{i} &=&\frac{1}{2}\Omega
_{jk}^{c}g_{cb}g^{ji}-q_{jk}^{ih}C_{hb}^{j},\ _{\shortmid }Z_{jb}^{a}=-\
^{\pm }q_{cb}^{ad}\ \Xi _{dj}^{c},\ _{\shortmid }Z_{bc}^{a}=0,  \notag \\
\ _{\shortmid }Z_{bk}^{a} &=&\ ^{+}q_{cd}^{ab}\Xi _{bk}^{c},\ \ _{\shortmid
}Z_{kb}^{i}=\frac{1}{2}\Omega _{jk}^{a}g_{cb}g^{ji}+C_{hb}^{j}q_{jk}^{ih},\
_{\shortmid }Z_{jk}^{i}=0,  \notag \\
\ _{\shortmid }Z_{ab}^{i} &=&-\frac{g^{ij}}{2}\left\{ g_{cb}\ \Xi
_{aj}^{c}+g_{ca}\ \Xi _{bj}^{c}\right\} ,\ _{\shortmid
}Z_{jk}^{a}=-C_{jb}^{i}g_{ik}g^{ab}-\frac{1}{2}\Omega _{jk}^{a},
\label{cdeftc}
\end{eqnarray}%
for $\ q_{jk}^{ih}=\frac{1}{2}(\delta _{j}^{i}\delta
_{k}^{h}-g_{jk}g^{ih}),\ ^{\pm }q_{cd}^{ab}=\frac{1}{2}(\delta
_{c}^{a}\delta _{d}^{b}\pm g_{cd}g^{ab}),\ \Xi _{aj}^{c}=\left[
L_{aj}^{c}-e_{a}(N_{j}^{c})\right] .$

If $\mathbf{V}=TM,$ for certain models of Lagrange and/or Finsler geometry,
we can identify $\widehat{L}_{jk}^{i}$ to $\widehat{L}_{bk}^{a}$ and $%
\widehat{C}_{jc}^{i}$ to $\ \widehat{C}_{bc}^{a}$ and consider the canonical
d--connection as a couple $\widehat{\mathbf{\Gamma }}_{\ \alpha \beta
}^{\gamma }=\left( \widehat{L}_{jk}^{i},\widehat{C}_{jk}^{i}\right) .$

By a straightforward d--form calculus, we can find the N--adapted components
of the curvature of a d--connection $\mathbf{D},$
\begin{equation}
\mathcal{R}_{~\beta }^{\alpha }\doteqdot \mathbf{D\Gamma }_{\ \beta
}^{\alpha }=d\mathbf{\Gamma }_{\ \beta }^{\alpha }-\mathbf{\Gamma }_{\ \beta
}^{\gamma }\wedge \mathbf{\Gamma }_{\ \gamma }^{\alpha }=\mathbf{R}_{\ \beta
\gamma \delta }^{\alpha }\mathbf{e}^{\gamma }\wedge \mathbf{e}^{\delta };
\label{curv}
\end{equation}
\begin{eqnarray}
R_{\ hjk}^{i} &=&e_{k}L_{\ hj}^{i}-e_{j}L_{\ hk}^{i}+L_{\ hj}^{m}L_{\
mk}^{i}-L_{\ hk}^{m}L_{\ mj}^{i}-C_{\ ha}^{i}\Omega _{\ kj}^{a},  \notag \\
R_{\ bjk}^{a} &=&e_{k}L_{\ bj}^{a}-e_{j}L_{\ bk}^{a}+L_{\ bj}^{c}L_{\
ck}^{a}-L_{\ bk}^{c}L_{\ cj}^{a}-C_{\ bc}^{a}\Omega _{\ kj}^{c},  \notag \\
R_{\ jka}^{i} &=&e_{a}L_{\ jk}^{i}-D_{k}C_{\ ja}^{i}+C_{\ jb}^{i}T_{\
ka}^{b},  \label{dcurv} \\
R_{\ bka}^{c} &=&e_{a}L_{\ bk}^{c}-D_{k}C_{\ ba}^{c}+C_{\ bd}^{c}T_{\
ka}^{c},  \notag \\
R_{\ jbc}^{i} &=&e_{c}C_{\ jb}^{i}-e_{b}C_{\ jc}^{i}+C_{\ jb}^{h}C_{\
hc}^{i}-C_{\ jc}^{h}C_{\ hb}^{i},  \notag \\
R_{\ bcd}^{a} &=&e_{d}C_{\ bc}^{a}-e_{c}C_{\ bd}^{a}+C_{\ bc}^{e}C_{\
ed}^{a}-C_{\ bd}^{e}C_{\ ec}^{a}.  \notag
\end{eqnarray}

Contracting respectively the components of (\ref{dcurv}), one proves that
the Ricci tensor $\mathbf{R}_{\alpha \beta }\doteqdot \mathbf{R}_{\ \alpha
\beta \tau }^{\tau }$ is characterized by  d--tensors,%
\begin{equation}
R_{ij}\doteqdot R_{\ ijk}^{k},\ \ R_{ia}\doteqdot -R_{\ ika}^{k},\
R_{ai}\doteqdot R_{\ aib}^{b},\ R_{ab}\doteqdot R_{\ abc}^{c}.
\label{dricci}
\end{equation}%
It should be noted that this tensor is not symmetric for arbitrary
d--connecti\-ons $\mathbf{D}.$ The scalar curvature of a d--connection is
\begin{equation}
\ ^{s}\mathbf{R}\doteqdot \mathbf{g}^{\alpha \beta }\mathbf{R}_{\alpha \beta
}=R+S,\ \ R=g^{ij}R_{ij},\ S=g^{ab}R_{ab},  \label{sdccurv}
\end{equation}%
defined by a sum the h-- and v--components of (\ref{dricci}) and d--metric (%
\ref{m1}).

The Einstein tensor is defined and computed in standard form
\begin{equation}
\mathbf{G}_{\alpha \beta }=\mathbf{R}_{\alpha \beta }-\frac{1}{2}\mathbf{g}%
_{\alpha \beta }\ ^{s}\mathbf{R}  \label{enstdt}
\end{equation}%
It should be noted that, in general, this Einstein tensor is different from
that defined for the Levi--Civita connection but for the canonical
d--connection and metric defined by a Lagrange model both such tensors are
derived from the same Lagrangian and metric structure. Finally, we note that
formulas (\ref{curv})--(\ref{enstdt}) are defined in the same form for
different classes of linear connection. For the canonical d--connection and
the Levi--Civita connection, we label such formulas with respective 'hats'
and left 'vertical lines', for instance, $\widehat{\mathbf{R}}_{\ \beta
\gamma \delta }^{\alpha }$ and $\ _{\shortmid }R_{\ \beta \gamma \delta
}^{\alpha },\widehat{\mathbf{R}}_{\alpha \beta }$ and $\ _{\shortmid
}R_{\alpha \beta },...$

\end{document}